\documentclass[12pt]{article}
\usepackage{amsfonts}
\usepackage{mathrsfs}
\usepackage{amsmath}
\usepackage{latexsym}
\usepackage{amssymb}
\usepackage{amsthm}
\usepackage{amscd}

\title{The radicals of (Braided ) Lie algebras }

\author{  \small Lingwei Guo, \small Shouchuan Zhang,   \small Junqin Li \\
\small $a$. Department  of Mathematics, Hunan University\\
\small   Changsha  410082, \ P.R. China \\
\small Emails: z9491@yahoo.com.cn
 }

\date {}
\begin{document}
\newtheorem{Proposition}{\quad Proposition}[section]
\newtheorem{Theorem}[Proposition]{\quad Theorem}
\newtheorem{Definition}[Proposition]{\quad Definition}
\newtheorem{Corollary}[Proposition]{\quad Corollary}
\newtheorem{Lemma}[Proposition]{\quad Lemma}
\newtheorem{Example}[Proposition]{\quad Example}
\newtheorem{Remark}[Proposition]{\quad Remark}

\maketitle \addtocounter{section}{-1}

\numberwithin{equation}{section}

\begin {abstract} The general theory of the radicals of Lie algebras are established. Baer radicals of untwisted affine Lie algebras are found.

\vskip0.1cm 2000 Mathematics Subject Classification: 17B65,  17B67.

keywords: Radical,  Lie algebra.
\end {abstract}

\section{Introduction}\label {s0}
In order to study infinite dimensional Lie algebras V. Kac and R.
Moody independently introduced  Kac-Moody algebras in 1960s ( \cite
{Ka85, Wa02, Xu07} ). The constructure of finite dimensional Lie
algebras are clear (\cite {Hu72, Me99, Wa78, Ja62}). However, we
know only a little for infinite dimensional Lie algebras (\cite
{OZ03, LT07, Ka85}). On the other wise, the radical theory of
associative rings is every complete (\cite {Di65, Sz82}). However,
it is almost the beginning for Lie algebras.

There exists a maximal solvable ideal for every finite dimensional
Lie algebra $L$, which is called the radical  of $L$. But it is
possible that the result above  does not hold for infinite
dimensional Lie algebras. However, we introduce the theory of
radicals of Lie algebras (contains infinite dimensional Lie
algebras) in this paper. The Baer radical exactly is   the maximal
solvable ideal for every finite dimensional Lie algebra.

In this paper, we built the radical theory of Lie algebras and
introduce Baer radical of Lie algebras. We obtain the Baer radicals
of Witt algebras, Virasoro algebra, loop algebras and  untwisted
affine algebras.

Let $F$ be a field of algebraic closure  with  characteristic zero
and ${\mathbb Z}$ the set of all integers.

\section{General theory of the radicals}\label {s0}
In this section, the general theory of the radicals of Lie algebras
is established.

If $\bar L$ is a homomorphic image of $L$, then we write $L \sim
\bar L$.
\begin{Definition}\label{1.1}
The property $r$ is called a radical property, if $r$ satisfies the following three axioms:

(R1) Every homomorphic image $\bar L$ of an $r$- Lie algebra is
again an $r$-Lie algebra;

(R2)Every Lie algebra $L$ has an $r$-ideal  $r(L)$ of $L$;

(R3)The factor Lie algebra $L/r(L)$ of $L$ with respect to $r(L)$ is $r$-semisimple, i.e. $r(L/r(L))=(0)$.
\end {Definition}

\begin{Lemma}\label{1.2} If $I, J \lhd L,$ then $(I+J)/ I \cong J/ I\cap J.$
\end{Lemma}
{\bf Proof.} Let $\varphi$ be a map from $I+J$ to $J/(I \cap J)$
sending $x+y$ to $y +I \cap J$ for all $x\in I, y\in J.$ It is clear
that  $\varphi$ is a Lie algebra homomorphism and $ker \varphi = I.$
$\Box$
\begin{Lemma}\label{1.3} If $I, J \lhd L,$  and $I \subseteq J$,  then $L/ J \cong (L/I)/(J/I).$
\end{Lemma}
{\bf Proof.} Let $\varphi$ be a map from $L/I$ to $L/J$ sending $x+
I$ to $x +J $ for all $x\in L.$ It is clear that $\varphi$ is a Lie
algebra homomorphism and $ker \varphi = J/I.$ $\Box$
\begin{Lemma}\label{1.4} If $r$ is a radical property and $r(L/I) =0$ with $I \lhd L$, then $r(L) \subseteq I.$
\end{Lemma}

{ \bf Proof. }Considering $(I+r(L))/I\cong r(L)/I\cap r(L)$, we have $(I+r(L))/I=0$,
which implies $r(L)\subseteq I$. $\Box$

\begin{Proposition}\label{1.5} If $r$ is a radical property and $\varphi $ is a Lie algebra
homomorphism from $L$ onto $L$, then $\varphi (r(L)) \subseteq
r(L)$.
\end{Proposition}
\begin{Proposition}\label{1.6}
If $B$ is an ideal of a Lie algebra $L$ and if $B$ and $L/B$ are
$r$-Lie algebras, then $L$ is also an $r$-Lie algebra.
\end{Proposition}
{\bf Proof.} Considering $B\subseteq r(L)$, we have that  $L/ r(L)
\cong (L/B)/ (r(L)/ B)$ is not only an $r$-ideal but also
$r$-semisimple. Consequently, $L/ r(L)=0$ and $L = r(L).$

\begin{Proposition}\label{1.7}
In a Lie algebra $L$ the union $A$ of any ascending chain of
$r$-ideals $\{L_i \mid i =1, 2, \cdots \}$ is again an $r$-ideal.
\end{Proposition}
{\bf Proof.} Since $(L_i + r(A))/ r(A) \cong L_i / (r (A) \cap
L_i)$, we have that $(L_i + r(A))/ r(A)$ is an $r$-ideal of $A/
r(A)$, which implies $L_i + r(A) = r(A)$ for any $i=1, 2, \cdots.$
Consequently, $r(A) =A.$ $\Box$

\begin{Proposition}\label{1.8}
If $A_\alpha $ is an $r$-ideal of $L$ for any $\alpha \in I$, then
$A := \sum _{\alpha \in I} A_\alpha $  is also an $r$-ideal of $L$.
\end{Proposition}
{\bf Proof.} For any $\alpha \in I$, $(A_\alpha + r(A))/ r(A) \cong
A_\alpha / (A_\alpha \cap r(A))$ is an $r$-ideal of $A/r(A)$.
Consequently, $(A_\alpha + r(A))/ r(A)=0$ and $A_\alpha \subseteq
r(A),$ which implies $A = r(A).$  $\Box$

The following follows from Proposition \ref {1.8}
\begin{Theorem}\label{1.9}
The $r$-radical $r(L)$ of an arbitrary Lie algebra $L$ is the sum of all $r$-ideals.
\end{Theorem}

\begin{Proposition}\label{1.10}
If $B,C$ and $D$ are ideals of a Lie algebra $L$ , such that
$B\supseteqq C$ and $B/C$ is an $r$-Lie algebra , then $(B+D)/(C+D)$
is also an $r$-Lie algebra.
\end{Proposition}
{\bf Proof.} $$ (B/C)/((B \cap (C+D))/C) \cong B/ (B \cap (C+D))
\cong (B +(C+D))/(C+D) = (B+D)/ (C+D). \ \ \Box$$

\begin{Proposition}\label{1.12}
If $B$ is an ideal of $L$, $r$ is a radical property, $B$ and $L/B$
are $r$-semisimple Lie algebras, then $L$ is also $r$-semisimple.
\end{Proposition}

{ \bf Proof. } By Lemma 1.4, $r(L)\subseteq B_\alpha$ for any $\alpha\in I$ and
$r(L)\subseteq B_0$. Set $r(L)\nsubseteq B_{\alpha_0}/B_0$. Set $r(L/B_0)=W/B_0$.
See $(W+B_{\alpha_0})/B_{\alpha_0}\cong W/(B_{\alpha_0}\cap W)$ and $B_0=W\cap B_0
\subseteq W\cap B_{\alpha_0}$. Since $(W+B_{\alpha_0})/B_{\alpha_0}\cong W/(B_{\alpha_0}\cap W)\cong (W/B_{\alpha_0})/((B_{\alpha_0}\cap W)/B_{\alpha_0})$,
we have that $(W+B_{\alpha_0})/B_{\alpha_0}$ is an $r$- ideal and $(W+B_{\alpha_0})/B_{\alpha_0}=0$, which implies a contradiction. $\Box$

\begin{Proposition}\label{1.14}
If $r$ is a radical property and $B_\alpha$ are ideals in the Lie
algebra $L$ such that each factor Lie algebra $L/B_\alpha$ is
$r$-semisimple and $B_0=\cap_{\alpha \in I}  B_\alpha$, then $L/B_0$
is also $r$-semisimple.
\end{Proposition}
{\bf Proof.} By Lemma \ref {1.4}, $r(L) \subseteq B_\alpha$ for any
$\alpha \in I$ and $r(L) \subseteq B_0.$ If $r(L/ B_0) \not= 0$,
then there exists $\alpha_0 \in I$ such that $r(L / B_0) \nsubseteq
B_{\alpha_0}/ B_0$. Set $r(L/ B_0) = W/B_0.$  Since $$(W +B_{\alpha
_0})/ B_{\alpha_0} \cong W/ (B_{\alpha_0} \cap W )\cong (W/ B_{0})/
(( B_{\alpha_0 } \cap W)/ B_{0})$$ we have that $0\not= (W
+B_{\alpha _0})/ B_{\alpha_0}$ is an $r$-ideal. This is a
contradiction. $\Box$
\begin{Theorem}\label{1.13}
The radical $r(L)$ of $L$ is the intersection $D$ of all ideals
$I_\alpha$ of $L$ such that $L/I_\alpha$ is $r$-semisimple. That is,
$r(L) = \cap \{ I \lhd L \mid r( L/I)=0 \}$.
\end{Theorem}

{ \bf Proof. } By Lemma 1.4, $r(L)\subseteq D$. On the other hand, $r(L)\lhd L$ and
$L/r(L)\in S(r)$, which implies $D\subseteq r(L)$. Therefore, $r(L)=D$. $\Box$

\begin{Theorem}\label{1.14}
The Lie algebra property $r$ is a radical property if and only if

$(R1')$ Every homomorphic image $\bar L$ of an $r$-Lie algebra $L$
is an $r$-Lie algebra.

$(R2')$ If every non-zero homomorphic image of a Lie algebra $L$ contains a non-zero $r$-ideal, then $L$ is an $r$-Lie algebra.
\end{Theorem}

{ \bf Proof. }Assume that r is a radical property. If every non-zero homomorphic image of
a Lie algebra L contains a non-zero $r$-ideal, then $r(L)=L$ since $L/r(L)$ has not any non-zero $r$-ideal. Consequently, $(R2')$ holds.

Assume that both $(R1')$ and $(R2')$ hold. Let $S:=\sum\{N \mbox{is
an r-ideal of } L\}$. If $S=0$, then S is an $r$-ideal of L. If
$S\neq0$ and $U\lhd S$ with $U\neq S$, then there exists an
$r$-ideal V of L such that $V\nsubseteq U$. Since $(V+U)/U\cong
V(U\cap V)$ and $V\cap U\neq V$, we have $(V+U)/U$ is a non-zero
$r$-ideal of $S/U$. By $(R2')$, S is an $r$-ideal. Consequently,
$(R2)$ holds.

If $r(L/r(L))\neq0$, then there exists $W\lhd L$ such that $r(L/r(L))=W/r(L)$
with $r(L)\neq W$. For any $X\lhd W$ with $W\neq X$, if $S\subseteq X$, then
$W/X\cong(W/S)/(X/S)$ is an $r$-ideal; if $S\nsubseteq X$, then $(X+S)/X\cong S/(X\cap S)$ is a non-zero $r$-ideal of $W/X$. Consequently, W is an $r$-ideal
of L by $(R2')$. This is a contradiction and $(R3)$ holds. $\Box$

\begin{Theorem}\label{1.15}
For a radical property $r$ the non-zero Lie algebra $L$ is an
$r$-Lie algebra if and only if $L$ can not be homomorphically mapped
onto a non-zero $r$-semisimple Lie algebra.
\end{Theorem}
{\bf Proof.} It follows from Theorem \ref {1.14}. $\Box$
\begin{Theorem}\label{1.16} Let ${\mathbb K}$ be a class of Lie algebras, and satisfies
 the following two conditions :

(Q1) Every non-zero ideal of every  Lie algebra in ${\mathbb K}$ can be homomorphically mapped onto non-zero Lie algebra in ${\mathbb K}$.

(Q2) Every non-zero ideal of every  Lie algebra $L$ can be homomorphically mapped  onto non-zero Lie algebra in ${\mathbb K}$, then $L\in {\mathbb K}.$\\
then there exists a radical property $r$ such that $S(r) = {\mathbb K}.$
\end{Theorem}
 {\bf Proof.}
 Define that Lie algebra $L$ is an $r$-Lie algebra if and only if $L$ can not  be homomorphically mapped  onto non-zero Lie algebra in ${\mathbb K}$. Now we show that $r$ is a radical property.

 (i) For ($R1'$), assume that  $L$ is an $r$- Lie algebra and $L\sim \bar L \not=0.$ If $\bar L$ is not $r$-Lie algebra, then $\bar L$ can be homomorphically mapped  onto non-zero Lie algebra in ${\mathbb K}$, which implies that $ L$ can be homomorphically mapped  onto non-zero Lie algebra in ${\mathbb K}$. This is a contradiction. Consequently, ($R1'$) holds.

 (ii) ${\mathbb K} = S(r).$ In fact, if $L \in {\mathbb K}$ and $L \notin S(r)$, then $0\not= r(L)$. By (Q1), $ r(L)$ can be homomorphically mapped  onto non-zero Lie algebra in ${\mathbb K}$, which implies a contradiction.   If $L\in S(r)$ and $L\notin {\mathbb K}$,  then there exists a non-zero ideal $B$ of $L$ such that $B$ can not be  homomorphically mapped onto non-zero Lie algebra in ${\mathbb K}$. Therefore, $B$ is an $r$-Lie algebra, which implies that $L$ is not $r$-semisimple. This is a contradiction.

 (iii)  Assume   that every non-zero homomorphic image of  $L$ has a non-zero $r$-ideal of $L$. Consequently,  $L$ can not  be homomorphically mapped  onto non-zero Lie algebra in ${\mathbb K}$. This implies that $L$ is an $r$-Lie algebra, that is, ($R2'$) holds. $\Box$

 \begin{Proposition}\label{1.17} If ${\mathbb K}$ satisfies $(Q1)$ and $\bar {\mathbb K} := \{ L \mid $ every non-zero ideal of $L$ can be homomorphically mapped onto a non-zero Lie algebra in ${\mathbb K}\}$, then $\bar {\mathbb K}$ satisfies (Q1) and (Q2).
\end{Proposition}

 {\bf Proof.}It is clear ${\mathbb K}\subseteq\bar{\mathbb K}$. If $L\in \bar{\mathbb K}$ and $0\neq I\lhd L$,
 then there exists $0\neq\bar{I}\in{\mathbb K}$ such that $I\sim\bar{I}$. Consequently, $\bar{\mathbb K}$
 satisfies $(Q1)$. For $(Q2)$, if every non-zero I ideal of L can be homomorphically mapped onto a non-zero
 Lie algebra $\bar{I}$ in $\bar{\mathbb K}$, which implies that there exists $0\neq\bar{\bar I}\in{\mathbb K}$
 such that $\bar{I}\sim\bar{\bar I}$. Consequently, $I\sim\bar{\bar I}$ and $L\in\bar{\mathbb K}$. That is,
 $\bar{\mathbb K}$ satisfies $(Q2)$. $\Box$

Let $r ^{\mathbb K}$  denote the radical property  determined by $\bar {\mathbb K}$ and $r ^{\mathbb K}$ is  called the upper radical of ${\mathbb K}$.

 \begin{Proposition}\label{1.18} If both $r$ and $r'$ are  radical properties, then $R(r) \subseteq R(r')$ if and only if $S(r') \subseteq S(r)$.
\end{Proposition}

 {\bf Proof.} If $R(r)\subseteq R(r')$ and $L\in S(r')$, then L has not any non-zero
 $r'$-ideal. Consequently, L has not any non-zero $r$-ideal and $L\in S(r)$. That is, $S(r')\subseteq S(r)$.

 Conversely, if $S(r')\subseteq S(r)$ and $L\in R(r)$, then L can not be homomorphically mapped onto a non-zero $r$-semisimple Lie algebra by Theorem 1.15,
 which implies that L can not be homomorphically mapped onto a non-zero $r'$-semisimple Lie algebra. Consequently, $L\in R(r')$. $\Box$

\begin{Proposition}\label{1.19} If ${\mathbb K}$ satisfies $(Q1)$ and $r$ is a radical property with ${\mathbb K} \subseteq S(r)$, then $r \le r^{\mathbb K}$.
\end{Proposition}

{\bf Proof.}First we show $\bar{\mathbb K}\subseteq S(r)$, where $\bar{\mathbb K}$
is defined in Proposition 1.17. If $L\in \bar{\mathbb K}$ and $r(L)\neq0$, then
there exists $0\neq I\in{\mathbb K}$ such that $r(L)\sim I\in S(r)$, which is a
contradiction. Consequently, $r(L)=0$. $\Box$

\begin{Definition}\label{1.20} $A$ is called an accessible ideal of Lie algebra $L$ if there exists $A_1 \lhd A_2 \lhd \cdots \lhd A_n = L$ such that $A_1 = A$, written $A \lhd \lhd L.$
\end{Definition}

\begin{Proposition}\label{1.21} If ${\mathbb K}$ is a class of Lie algebras and
$\bar {\mathbb K} := \{ L \mid  L$ is an accessible ideal of a Lie algebra in  ${\mathbb K}\}$, then $\bar {\mathbb K}$ satisfies (Q1).
\end{Proposition}

\section{ The  special radicals for Lie algebras}\label {s2}

  In this section,
  special radicals for Lie algebras is defined.

Lie algebra $L$ is called a simple Lie  algebra if $L$ has not any
non-trivial  ideals and $L^2 \not=0.$
 $L$ is said to be  semiprime if there are no
non-zero nilpotent  ideals in $L$.  $L$ is said to be
  prime
if $[I, J] =0$ implies $I=0$ or $J=0$ for any ideals $I$ and $J$ of
$L$. A ideal $I$ is called an   (semi)prime ideal of $L$ if $L/I$ is
(semi)prime.

\begin{Lemma}\label{2.1}
   If  $L$  is an  Lie algebra and $E$ is a non-empty subset
of $L$, then
   $(E) = \sum _{k=0}^\infty ({\rm ad } L)^k E$,
    where $(E)$ denotes the  ideal generated by $E$ in  $L$.
\end {Lemma}
{ \bf Proof. } It is clear. $\Box$
  \begin{Proposition}\label{2.2}
  \label {9.1.2}   (i)  $L$ is  semiprime iff $( a)^2 = 0$ always implies
    $a = 0$ for any $a\in L$.

 (ii) $L$ is  prime iff $ [(a), (b)] = 0$ always implies
 $a = 0$  or $b = 0$ for any $a, b \in L$.
\end{Proposition}
{\bf Proof.}
It is clear. $\Box$

 \begin{Proposition}\label {2.3}
  If $ I \lhd L$ and $I$ is an  semiprime Lie algebra,
  then
       $I \cap I^{*} = 0$ and   $I^{*} \lhd L$,
     where $I^* = \{ a \in L \mid  [I, ( a)] = 0 \}$.

 \end{Proposition}
 {\bf Proof .} For any $x \in I^* \cap I$, we have that  $[I, ( x)] = 0$
 and  $[(x), ( x)] = 0$. Since $I$ is an  semiprime
Lie  algebra,  $x =0$,
 i.e. $I \cap I^* =0$.

 For any $x \in I^*, z \in L,$  see that
\begin {eqnarray*}[([z, x] ), I] &\subseteq & \sum _{k=0}^\infty [(ad L)^k [z, x],  I]\\
&\subseteq& [(x), I] =0.\end {eqnarray*} Therefore, $[z, x] \in I^*.$ That is, $I^*$ is an ideal of $L.$ $\Box$

 \begin{Definition}\label {2.4'} ${\mathbb K}$ is called a  weakly  special class if

     (WS1)  ${\mathbb K}$ consists of  semiprime Lie algebras .

     (WS2)  For any $L \in {\mathbb K}$, if $0 \not= I \lhd L$  then
     $I \in {\mathbb K}$.

     (WS3)  If $L$ is a  Lie  algebra and  $B\lhd L$ with
      $B \in {\mathbb K}$ and $B^* =0$,
  then    $ L \in{\mathbb  K}$,
where  $B^{*} = \{ a \in L\mid  [( a), B] = 0 \}$.
    \end{Definition}

 \begin{Definition}\label {2.4} ${\mathbb K}$ is called a special class if

     (S1)  ${\mathbb K}$ consists of  prime Lie algebras .

     (S2)  For any $L \in {\mathbb K}$, if $0 \not= I \lhd L$  then
     $I \in {\mathbb K}$.

     (S3)  If $L$ is a  Lie  algebra and  $B\lhd L$ with
      $B \in {\mathbb K}$,
  then    $ L/B^{*}\in{\mathbb  K}$.
    \end{Definition}

    By (WS2) or (S2), ${\mathbb K}$ satisfies $(Q1).$ It is clear
    that every special class is a weakly special class.

    Example. ${\mathbb K} := \{ L \mid L \hbox { is a simple Lie algebra}\}
    $ is a special class.

\section{ Baer radical}\label {s3}

In this section, we give the characterization of  Baer
radical.

Let   $(I, <)$ be  a set with well order. If $ \alpha \in I$ and
there exists $\gamma \in I$ such that $\beta  < \gamma < \alpha$ for
any $\beta \in I$ with $\beta < \alpha$, then $\alpha $ is called a
limit number. If $ \alpha \in I$ and  there  exists a $\beta \in I$
with $\beta < \alpha$ such that there does not exist any elements in
between $\alpha$ and $\beta$, then  $\alpha $ is called a non limit
number (or say that $\alpha$ is not a limit number). In this case,
$\beta = \alpha -1$, or $\alpha = \beta +1$.

Let ${\mathbb E}$ be a non-empty class of Lie algebras. Let $I$ be a well order set with the maximal element $ \alpha_0$ of $I$. Define ${\mathbb E}_1 :=\{L \mid \hbox { there exists } \bar L \in {\mathbb E}   \hbox{ such that } \bar L \sim L \}  $.
 Assume that  for any $\beta < \alpha $ with $\alpha , \beta \in I$,
 ${\mathbb E}_\beta  $ is defined. If $\alpha $ is a limit  number,
 define $E_\alpha := \{L \mid \hbox { there exists  } \beta \in
 I  \hbox {  with }   \beta < \alpha
 \hbox{ such that } L \in  {\mathbb E}_\beta  \}  $; If $\alpha $ is not
  a limit  number, define $E_\alpha := \{L \mid  \hbox { every non-zero  homomorphic image
   of  } L   \hbox{ has a non-zero ideal }  \in  {\mathbb E}_{\alpha-1}  \}.$
  Let ${\mathbb E}_I: = \cup _{\alpha \in I} {\mathbb E}_\alpha $.

 \begin{Lemma}\label {3.1} For any $\alpha, \beta \in I$ with $\alpha \le \beta$,

 (i) ${\mathbb E}_\alpha$ is closed under homomorphisms.

 (ii) ${\mathbb E}_\alpha \subseteq {\mathbb E}_\beta.$

   \end{Lemma}

{\bf Proof.} (i) Assume  $L \in {\mathbb E}_\alpha$ and $L \sim \bar L$. We show $\bar L \in {\mathbb E}_\alpha$ by transfinite  induction.
If $\alpha =1,$ then there exists $\widetilde{ L} \in {\mathbb E}$ such that $\widetilde{L} \sim L$. Consequently, $\widetilde{L} \sim \bar L$, which implies $\bar L
\in {\mathbb E}_1.$ Assume $\alpha >1$. If $\alpha $ is a limit number, then there exists $\beta \in I$ such that $\beta < \alpha$ and $L \in {\mathbb E}_\beta$, which implies
$\bar L \in {\mathbb E}_\alpha$. If $\alpha $ is not  a limit number, then every non-zero homomorphic image of $L$ has a non-zero ideal in ${\mathbb E}_{\alpha-1}$. Consequently, every non-zero homomorphic image of $\bar L$ has a non-zero ideal in ${\mathbb E}_{\alpha-1}$, which implies  $\bar L \in {\mathbb E}_\alpha.$

(ii) By Part  (i), ${\mathbb E}_\gamma \subseteq {\mathbb
E}_{\gamma+1}$ for any $\gamma, \gamma +1 \in I$. If $\beta $ is a
limit number, then the claim holds. Now we assume that  $\beta$ is
not a limit number and  show  (ii)  by transfinite induction on
$\alpha$. For $\alpha =1 $ and $\beta
>1$,  By
inductive assumption on $\beta$, ${\mathbb E}_1 \subseteq {\mathbb
E}_{\beta -1} $ since $1\le \beta -1 $. Consequently, ${\mathbb E}_1
\subseteq {\mathbb E}_{\beta } .$ Now assume $\alpha
>1.$ If   $L \in {\mathbb
E}_\alpha$, then every non-zero homomorphic image of $ L$ has a
non-zero ideal in ${\mathbb E}_{\alpha-1}\subseteq {\mathbb
E}_{\beta -1} $. Therefore, ${\mathbb E}_{\alpha}\subseteq {\mathbb
E}_{\beta } $.   $\Box$

 \begin{Theorem}\label {3.2} Let $I$ and $I'$ be two well order sets with the maximal elements. If ${\rm card} I > {\rm card} 2^L$ and  ${\rm card} I' > {\rm card} 2^L$, then
$L \in {\mathbb E_I}$ if and only if $L \in {\mathbb E_{I'}}$.

   \end{Theorem}
   {\bf Proof.} Let $\beta_0$ be the maximal of $I'$. Assume that $L \in {\mathbb E}_I$. Then there exists an $ \alpha \in I$ such that $L\in {\mathbb E}_\alpha$. Now we show that  there exists  $ \alpha' \in I'$ such that $L\in {\mathbb E}_{\alpha '}$.
 If $\alpha =1$, then the conclusion holds. Assume $\alpha >1.$ If $\alpha$ is a limit number, then there exists $ \beta \in I$ with $\alpha > \beta$ such that $L \in {\mathbb E}_\beta$. Therefore, the conclusion holds by induction assumption. If $\alpha $ is not a limit number,  then every non-zero homomorphic image $\bar L$  of   $L$    has a non-zero ideal $ N  \in  {\mathbb E}_{\alpha-1}$.
 By induction assumption, there exists $\gamma \in I'$ such that $N \in {\mathbb E}_\gamma.$ Let $\{ N_\mu \mid  \mu \in J  \}$ is the set of all ideals of $L$
   with $J \subseteq I'\setminus \beta_0$ and $0\not= M_\mu \lhd L/N_{\mu} $ with $M_\mu \in   {\mathbb E}_\mu$. If $\beta_0$ is a limit number,  $L \in  {\mathbb E}_{\beta_0}$. If $\beta_0$ is not a limit number, then $M_\mu \in   {\mathbb E}_{\beta_0 -1}$ for all
   $\mu \in J.$
  Therefore,  $L \in  {\mathbb E}_{\beta_0}.$ $\Box$

 $L$ is called an $r_{\mathbb E}$-Lie algebra if and only if there exists a  well order set with
 ${\rm card} I > {\rm card} 2^L$ and  the maximal $\alpha_0$  such that
 $L \in {\mathbb E}_I.$

 \begin{Theorem}\label {3.3}
$r_{\mathbb E}$ is a radical property.
   \end{Theorem}

  {\bf Proof.} Let $r = r_{\mathbb E}$ for convenience.

   (R1') If $L$ is an $r$-Lie algebra and  $L \sim \bar  L$, then  there exists $\alpha \in I$ such that $L \in {\mathbb E}_{\alpha}$. By Lemma \ref {3.1}, $\bar L \in {\mathbb E}_{\alpha}.$

   (R2') Assume that every non-zero homomorphic image of  $L$ has a non-zero $r$-ideal.
  Let $I$ be a  well order set with
 ${\rm card} I > {\rm card} 2^L$ and  the maximal $\alpha_0 $ and $\alpha _0 -1 \in I$.
 Let $J \subseteq I$ and $\alpha _0, \alpha _0 -1 \notin J$  such that $\{ N_i \in i \in J\}$ is the set of all ideals of $L$. Consequently, there exists $0\not= B_i \lhd L/ N_i$  such that  $B_i $ is $r$-ideal of  $L/ N_i$ for any $i\in J$. Thus $B_i \in {\mathbb E}_{\alpha_0-1}$ for any $i \in J$ by Lemma \ref {3.1} since $i < \alpha _0 -1$ for any $i\in J$. Consequently, $L \in {\mathbb E}_{\alpha_0 }$, which implies
 $L \in {\mathbb E}_I.$ $\Box$

 If ${\mathbb E} = \{L \mid L \hbox { is a nilpotent Lie algebra}\}$, then $r_{\mathbb E}$ is called the Baer radical, written as $r_b$.

 \begin{Proposition}\label {3.4}
 (i) If $r(L)=0$, then $L$ has not any non-zero nilpotent ideals.

 (ii) If $L$
has not any accessible non-zero nilpotent ideals, then $r_b (L)=0.$

(iii) If $L$ is a solvable Lie algebra, then $L$ is an $r_b$-Lie
algebra.
   \end{Proposition}
   {\bf Proof.} (i) It is clear.

   (ii) If $r_b(L) \not=0,$ then there exists $\alpha \in I$ such that $r_b(L) \in {\mathbb E}_\alpha.$ We show that there exists an accessible non-zero nilpotent ideal
   of $L$ by transfinite induction on $\alpha.$ If $\alpha =1,$ then $r_b(L)$ is a non-zero nilpotent ideal of $L$. Assume $\alpha >1.$ If $\alpha $ is a limit number, then the claim holds. If $\alpha $ is not a limit number, then there exists $0\not= B \lhd r_b(L)$ such that $B\in {\mathbb E}_{\alpha -1}$. Consequently, $B$ has an accessible non-zero nilpotent ideal, which is also  an accessible non-zero nilpotent ideal of $L$.

(iii) On one hand, $L/r_b (L)$ is $r_b$-semisimple. On the other
hand, $L/r_b (L)$ is solvable. Consequently, $L/r_b (L)=0. \ \Box$

 \begin{Definition}\label {3.5} Let  $L$ be a Lie algebra and $I$ a well order set with ${\rm card} I > {\rm card} 2^L$, which has the maximal number $\alpha_0$ of $I$ and is not limit number. Define $L_{1} =  \sum \{ A \mid A \lhd L, A \in {\mathbb E}_1 \}$.
 Assume $\alpha >1$. If $\alpha $ is a limit number, define $L_\alpha :=
 \sum \{ A \mid A \lhd L, \hbox { there exists } \beta \in I, \beta < \alpha, \hbox  { such that } A \in   {\mathbb E} _\beta \}$. If  $\alpha $ is not  a limit number, define $L_\alpha/ L_{\alpha -1} :=
 \sum \{ A / L_{\alpha -1} \mid A \lhd L,  A/ L_{\alpha -1} \in   {\mathbb E} _{\alpha -1} \}$.  Let $\{N_i \mid i \in J\}$ be the set of all ideal of $L$ with $J \subseteq I$ and $\alpha _0, \alpha _0-1\notin J$,
 since  ${\rm card} I > {\rm card} 2^L$. Thus it is clear
 $L_{\alpha_0 -1} = L_{\alpha _0}$. Consequently,    there exists $\tau \in I$ such that $L_\tau = L_{\tau +1} = L_{\tau +2}= \cdots.$ Written $L_\tau = L_{\mathbb E}.$
 \end{Definition}

   \begin{Proposition}\label {3.6} $r_{\mathbb E} (L) = L_{\mathbb E}$.
    \end{Proposition}
 {\bf Proof.}  Let $L_{\mathbb E} = L_\tau$.

 We first show taht  $L _\alpha \subseteq  r _{\mathbb E} (L)$ by transfinite
 induction on $\alpha$. It is clear $L_1 =
  \sum \{ I \lhd L \mid I \in {\mathbb E}_1\}\subseteq  r _{\mathbb E}
  (L)$. Now assume $\alpha >1.$ If $\alpha$ is a limit number, then
  $L_\alpha :=
 \sum \{ A \mid A \lhd L, \hbox { there exists } \beta \in I,
 \beta < \alpha, \hbox  { such that } A \in   {\mathbb E} _\beta \}\subseteq r _{\mathbb E}
  (L)$ by inductive assumption. If   $\alpha $ is not  a limit number,
  then \begin {eqnarray*}L_\alpha/ L_{\alpha -1} &=&
 \sum \{ A / L_{\alpha -1} \mid A \lhd L,  A/ L_{\alpha -1} \in
  {\mathbb E} _{\alpha -1} \} \subseteq r_{\mathbb E} (L/ L_{\alpha -1}). \end
  {eqnarray*} Therefore $L_\alpha \subseteq r_{\mathbb E} (L)$ by Proposition \ref {1.6}.

  Next we show $r_{\mathbb E}(L) \subseteq  L_{\tau}$.
 Considering Lemma \ref {1.4},  we only need
  show $r_{\mathbb E}(L/ L_{\tau})=0$. If $r_{\mathbb E}(L/ L_{\tau})\not=0$,
  then there exists $\alpha \in I$ and $A \lhd L$ such that $r_{\mathbb E} (L/L_\tau)=A/L_{\tau}
  \in {\mathbb E}_\alpha$. If $\alpha \ge \tau$, then  $A/ L_{\alpha} =
   A/ L_{\tau} \in {\mathbb E}_\alpha.$  Thus $A \subseteq L_{\alpha +1} $,
   i.e. $L_{\alpha +1} \not= L_\alpha$. This is a contradiction. If  $\alpha < \tau$, then $ A/ L_{\tau} \in {\mathbb E}_\alpha \subseteq {\mathbb E}_\tau.$ This is a contradiction.
 $\Box$

\section {Levi decomposition}\label {s4}
In this section we show that if $L$ is a finite dimensional Lie
algebra, then
  $r(L)$ is the maximal solvable  ideal of $L$ and $L = S \oplus r_b (L),$ where $S$ is a semi-simple Lie algebra.

  \begin{Proposition}\label {4.1}  If $L$ is a simple Lie algebra and $r$ is a radical
  property, then $r_b(L) =0.$
\end{Proposition}

 \begin{Proposition}\label {4.2}  If $L = \oplus _{i\in I} L_i$ is a directed sum as Lie algebras
 and $r$ is a radical
  property, then $r(L)= \oplus _{i\in I} r(L_i)$.
\end{Proposition}
{\bf Proof.} Since $r(L_i) \lhd L$, we have $r(L_i) \subseteq r(L) $
for any $i\in I.$  Conversely, since $\pi _i (r(L)) \subseteq
r(L_i)$ for all $i\in I$, we have $r(L) \subseteq \oplus _{i\in I}
r(L_i)$. $\Box$

 \begin{Proposition}\label {4.3}  If $L$ is a finite dimensional Lie algebra, then
  $r_b(L)$ is the maximal solvable  ideal of $L$ and
 $L = S \oplus r_b (L)$ as vector spaces,   where $S$ is a semi-simple Lie algebra.
\end{Proposition}
{\bf Proof.} By \cite [P 91]{Ja62}, $L = S \oplus R$ as vector
spaces, where $S$ is a semi-simple Lie algebra and $R$ is the
maximal solvable ideal of $L$. Since $r_b (L/R) =0$, we have $r_b(L)
\subseteq R$ by Lemma \ref {1.4}. Obviously, $R \subseteq r_b(L)$.
Thus $R = r_b(L)$. $\Box$

\section {Semi-directed sum} \label {s5}
In this section we show that  if $r$ is a radical property and $r(L)=0$, then $L$ is a semi-directed sum of $\{ L/I \mid I \lhd L, r(L/I)=0\}$.

Let $L:=\prod _{\alpha\in W}L_\alpha$ be the directed product of
$\{L_\alpha \mid \alpha \in W\}$. Let $\pi_\alpha$ and
$\iota_\alpha$ denote a canonical projection and a canonical
injection of the directed product $\Pi_{\alpha\in W}L_\alpha$. If
$B$ is a Lie sub-algebra of $L$ and $\pi_\alpha (B) = L_\alpha$ for
any $\alpha \in W,$ then $B$ is called a semi-directed sum of
$\{L_\alpha \mid \alpha \in W\}$.

\begin{Proposition}\label {5.1}  If  $B$ is  a semi-directed sum of $\{L_\alpha \mid \alpha \in W\}$, then  $\cap _{\alpha \in W} I_\alpha =0$ and  $B/ I_\alpha\cong L_\alpha$ with $I_\alpha = ker \pi _\alpha$ for any $\alpha \in W.$
\end{Proposition}
{\bf Proof.} It is clear. $\Box$
\begin{Proposition}\label {5.2}  If $B$ is a Lie algebra with $I_\alpha \lhd B$
for any $\alpha \in W$ and $\cap _{\alpha \in W} I_\alpha=0$, then
 $B$ is homomorphic to  a semi-directed sum  of  $ \{B/I_\alpha \mid \alpha \in
 W\}$ as Lie algebras.
\end{Proposition}
{\bf Proof.} For any $\alpha \in W,$ Let $\varphi_\alpha$ be the
 canonical homomorphism  from $B$ to $A_\alpha := B/I_\alpha$ and $\varphi := \prod _{\alpha
 \in W} \varphi _\alpha$   a homomorphism from $B$ to $\prod _{\alpha \in W}
 A_\alpha$. Since $ker \varphi = \cap _{\alpha \in W} ker \varphi _\alpha
 =0$, we have $B \cong \varphi (B)$ and $\varphi (B)$ is a
 semi-directed sum of $ \{B/I_\alpha \mid \alpha \in
 W\}$. $\Box$

\begin{Proposition}\label {5.3}  If $r$ is a radical property and $r(L)=0$, then $L$ is a semi-directed sum of $\{ L/I \mid I \lhd L, r(L/I)=0\}$.
\end{Proposition}
{\bf Proof.} By Theorem \ref {1.15},  $r(L) = \cap \{I \mid I \lhd L, r(L/I)=0\}$. Consequently, our claim  holds by Proposition \ref {5.2}. $\Box$

\section {Examples} \label {s6}
In this section  Baer radicals of untwisted affine Lie algebras are found.

\subsection { Central extension }

Assume that $(L, [ , ]_0)$ is a Lie algebra and $\psi$ is a
2-cocycle on $L$. If $\bar L := L \oplus Fc$ as vector spaces and $[
a+\lambda c , b+ \lambda' c] := [a, b]_0 + \psi (a, b)c$ for any $a,
b \in L, \lambda , \lambda ' \in F, $ then $\bar L$ is called a
central extension of $L$.

 \begin{Lemma}\label {6.0}
 If $\bar L = L + Fc$ is a central extension of $L$ and $r_b(L) =0$, then
 $r_b(\bar L) = Fc.$  \end{Lemma}

 {\bf Proof.} It is clear that $\pi$ is a Lie algebra homomorphism
 from $\bar L$ to $L$ by sending $a +\lambda c$ to $a$ for any $a\in L, \lambda \in
 F.$ Consequently, $\bar L / Fc \cong L$ and $r_b (\bar L/ Fc) =0$,
 which implies $r_b (\bar L) \subseteq Fc.$ It is clear $Fc \subseteq r_b (\bar
 L).$ $\Box$

\subsection { Witt algebras and Virasoro algebras}
 Let $ F[t, t^{-1}]$ be a  Laurent polynomial and  $d_i :=
t^{i+1}\frac {d}{dt}$
 for all $i\in \mathbb {Z}.$ Obviously,  \begin {eqnarray*} \label {e1}[d_i, d_j]
= (j-i) d_{i+j}.\end {eqnarray*}
  \begin{Definition}\label {6.1} (see \cite {Xu07}) If $L := \sum _{i=-1}^\infty Fd_i$, then $L$ is called
   Witt algebra on $F[t]$. If $L := \sum _{i=-\infty}^\infty Fd_i$,
   then $L$ is called  Witt algebra on $F[t, t^{-1}]$. Both of the two lie algebras are called Witt algebras.
\end{Definition}

 \begin{Proposition}\label {6.2} Witt algebras are simple algebras and $r_b$-semisimple. \end{Proposition}
{\bf Proof.} (i) Let $L$ is  Witt algebra on $F[t]$. Assume $0\not=
I \lhd L$ and $0 \not= u = \sum _{i=-1}^m k_i d_i  \in I$ with $k_m
\not=0$ and $m$ is the minimal. If $m \not= -1$, then
\begin {eqnarray*}
[u, d _{-1}] &=& \sum _{i=-1}^m k_i [d_i, d_{-1}]\\
&=& \sum _{i=0}^m k_i (-1-i) d_{i-1} \in I,
\end {eqnarray*} which contradicts to that $m$ is the minimal.
Consequently,  $m =-1$ and $ d_{-1} \in I$,which implies that $I
=L.$

(ii) Let  $L$ be   Witt algebra on $F[t, t^{-1}]$. Assume $0\not= I \lhd L$.
It is clear that there exists  $0 \not= u = \sum _{i=-1}^m k_i d_i
\in I$ with $k_m \not=0$ and $m$ is the minimal. We obtain $d_{-1}
\in I$ similar to the proof of Part (i). This can complete the
proof. $\Box$

 \begin{Definition}\label {6.2'} (see \cite {Xu07})
Let  $(L, [,]_0)$  be  Witt algebra on $F[t, t^{-1}]$ and  $\bar L$
a central extension of $L$ as follows:
$$ [d_i +\lambda c, d_j + \lambda' c] = [d_i, d_j]_0 + \frac {i ^3 -i} {12} \delta _{i+j, 0}\ c.$$
$\bar L$ is called Virasoro algebra. \end{Definition}

 \begin{Corollary}\label {6.2''} Let $\bar L= L + Fc$ be Virasoro algebra. then $r_b(\bar L)
 = Fc$. \end{Corollary}
{\bf Proof.} It follows from Lemma \ref {6.0} and Proposition \ref
{6.2}. $\Box$
 \subsection { Loop algebras}

 \begin{Definition}\label {6.3}(see \cite {Wa02})
$L(\stackrel {\circ} {L}):=F[t, t^{-1}]\otimes \stackrel {\circ}
{L}$ is called the loop algebra of finite dimensional simple Lie
algebra $\stackrel {\circ} {L},$ where $[t^m\otimes x,t^n\otimes
y]=t^{m+n}\otimes[x,y]$, for any $ m,n\in Z,x,y\in\stackrel {\circ}
{L}.$
\end{Definition}

 \begin{Lemma}\label {6.5}
 If $0\not= I \lhd \lhd L(\stackrel {\circ} {L})$, there exists $B \lhd F[t, t^{-1}]$ such that $0\not=  B \otimes \stackrel {\circ} {L} \subseteq I$.\end{Lemma}
{\bf Proof.} Let $I= I_1\lhd I_2\lhd \cdots \lhd I_n = L.$ We show
that there exists $0\not=B_i \lhd F[t, t^{-1}]$ such that $B_i
\otimes \stackrel {\circ} {L} \subseteq I_i$ for $1\le i \le n$ by
induction.  It is clear for   $n=1$.  Now assume $n>1.$ By induction
assumption, there exists an ideal $0\not= B_{2}$ such that $B_{2}
\otimes \stackrel {\circ} {L} \subseteq I_{2}$. Let $ 0\not= f(t)
\in B_{2}.$

(i) We first show that there exists $0\not= u = g(t)\otimes x \in
I.$ Let $\Phi$ be a root system of $\stackrel {\circ} {L}$ and
$\{e_\alpha, e_{-\alpha} \mid \alpha \in \Phi^+\} \cup \{H_1, H_2,
\cdots, H_n\}$ a Weyl basis of $\stackrel {\circ} {L}$. For $ u=
\sum _{\alpha \in \Phi ^+} (u_\alpha \otimes e _\alpha + u_{-\alpha}
\otimes e _{-\alpha} ) + \sum _{i=1}^n u_i' \otimes H_i \in I$, let
$l(u)$ denote the number of non-zero elements in $\{v_\alpha,
v_{-\alpha}, v_i' \mid \alpha \in \Phi^+, i=1, 2, \cdots, n\}$ and
$v$  an element in $I$ such that $l(v):= {\rm min } \{ l(u) \mid
0\not= u \in I\}$. If $l(v)=1$, the claim holds. Now assume  $l(v)>1
$.

(a) Assume   $v_{\alpha} =0$ for any $\alpha \in \Phi.$ Let $v =
\sum _{i=1} ^m w_i \otimes H_i'$, where  $w_1, w_2, \cdots, w_m$ is
linearly independent and $0\not=H_i' \in \eta$, Cartan subalgebra of
$\stackrel {\circ} {L}$. Thus there exists $\alpha \in \Phi$ such
that $\alpha (H_1') \not=0,$ which implies that $[ v, f(t) \otimes
E_\alpha ] = (\sum _{i=1} ^m f(t) w_i \alpha (H_i')) \otimes
E_\alpha \not=0$ and $l([ v, f(t) \otimes E_\alpha ])=1$ and $[ v,
f(t) \otimes E_\alpha ]\in I$. This is a contradiction.

(b)  If there exists $1\le i_0 \le n$  such that  $0\not= v_{i_0}' $
and $v_j =0$ for any $1\le j \le n$,  $i_0 \not= j,$ then there
exists $\alpha_0 \in \Phi$ such that $v_{\alpha _0} \not=0.$  Thus
$0 \not= l([v, f(t) \otimes H_{\alpha_0}])< l(v)$ and $[v, f(t)
\otimes H_{\alpha_0}]\in I$. This is a contradiction.

(c) If $v'_i =0$ for any $1\le i \le n$, then there exists $\alpha
_0\in \Phi $ such that $v_{\alpha _0} \not=0$.  $0\not= [v, f(t)
\otimes E_{-\alpha_0}] \in I$ is
 case  case (b), or $l([v, f(t) \otimes e_{-\alpha_0}]) =1.$

(ii) We show that $B \otimes \stackrel {\circ} {L} \subseteq I$,
where $B$ is an ideal generated by $g(t)$  of $F[t, t^{-1}]$. It is
clear that $\sum _{i=0} ^\infty (ad \stackrel {\circ} {L})^i x$ is
the ideal generated by $x$ of $\stackrel {\circ} {L}$. However
$\stackrel {\circ} {L} =\sum _{i=0} ^\infty (ad \stackrel {\circ}
{L})^i x  $ since $ \stackrel {\circ} {L}$ is a simple.  We show
that $F[t, t^{-1}] ^m g(t) \otimes (ad \stackrel {\circ} {L})^m x
\subseteq I $ by induction on $m$. It is clear for $m=0.$ Now $m>0.$
For any $u = [w, v] \in
 (ad \stackrel {\circ} {L})^m x$ and $h(t)= p(t)g(t) \in B$ with $w \in
 \stackrel {\circ}
{ L}$ and $v \in (ad \stackrel {\circ} {L})^{m-1} x $,  we have that
$h(t) \otimes u =
 [p(t) \otimes w, g(t) \otimes v] \in I$ since $g(t) \otimes v \in I$ by
 inductive assumption. Therefore this claim holds.
$\Box$

 \begin{Proposition}\label {6.6} $r_b (L(\stackrel {\circ} {L}))=0.$ \end{Proposition}
{\bf Proof.} By Proposition \ref {3.4} (ii), it is enough to show
that $L(\stackrel {\circ} {L})$ has not any non-zero nilpotent
accessible ideals. If $I$ is a non-zero nilpotent accessible ideal
of  $L(\stackrel {\circ} {L})$, then there exists $0\not= B \lhd
F[t, t^{-1}]$ such that $B\otimes \stackrel {\circ} {L} \subseteq I$
by Lemma \ref {6.5}. Let $0\not= g(t) \in B$. It is clear
$0\not=g(t)^m \otimes (\stackrel {\circ} {L}) ^m \subseteq (B\otimes
\stackrel {\circ} {L})^m  \subseteq   I^ m$. This is a
contradiction. $\Box$

\subsection { Untwisted Affine  algebras}
Let  $d:=t\frac{d}{dt}$ and $(t^m\otimes x, t^n \otimes y) _0 := \delta _{m+n, 0}(x, y)$
 for any $x, y \in \stackrel {\circ} {L}.$

\begin{Definition}\label {6.7} Let
$\widetilde{L}(\stackrel {\circ} {L}):=L(\stackrel {\circ} {L})+Fc $ and
 $[a+\lambda c, b+\mu c] := [a,b] + \psi (a, b)c$   for any $a, b \in L(\stackrel {\circ}{L})$, $\lambda, \mu \in F,$ where  $\psi (a, b) := (da, b)_0$ for any $a, b \in L(\stackrel {\circ} {L})$.
\end{Definition}

 \begin{Proposition}\label {6.8'}
$r_b( \widetilde{L}(\stackrel {\circ} {L}) = Fc$. \end{Proposition}
 {\bf
Proof.} It follows from Lemma \ref {6.0} and Proposition \ref {6.6}.
$\Box$

 \begin{Lemma}\label {6.8}
 For any $0\not=I \lhd \widetilde{L}(\stackrel {\circ} {L})$, there exists
 $I_1 \subseteq I \cap  {L}(\stackrel {\circ} {L})$ such that  $ I_1 \lhd
  L( \stackrel {\circ} {L})$, $Fc \subseteq I$ and $I = I_1 +Fc$.\end{Lemma}
{\bf Proof.}   There exists $0\not= u = u_1 + kc \in I$ with $u_1\in
\stackrel {\circ} {L}$
 and $k\in F$. If $u_1=0$, then  there exists $k\not=0$ and $c \in I,$ which implies $Fc \subseteq I.$ If $ u_1\not=0,$ then there exists $0\not=f(t)=\sum _{i=0}^m k_it^i \in F[t, t^{-1}]$ with $k_m\not=0$ and $0\not=x\in \stackrel {\circ} {L}$ such that $u_1 = f(t) \otimes x$ by means of the method similar  proof in Lemma \ref {6.5}. Consequently,
 $0\not=[u, t^{-m}\otimes x] = m(x, x)k_m \in I$, which implies $c\in I.$
  Let $I_1 := I \cap  L( \stackrel {\circ} {L}).$ It is clear that $ I_1 \lhd
  L( \stackrel {\circ} {L})$ and  $I = I_1+ Fc.$  $\Box$

 \begin{Definition}\label {6.9} (see \cite {Wa02})
$\widehat {L} (\stackrel {\circ} {L}):=\widetilde {L} (\stackrel
{\circ} {L}) +Fd$ is called untwisted affine Lie algebras, where $
[a + \lambda c + \lambda _1 d, a + \mu c + \mu _1 d]= [a, b]+
\lambda_1 db - \mu _1 da + \psi (a, b)c,$ \ \ where
$d=t\frac{d}{dt}$; $a, b \in {L} (\stackrel {\circ} {L}),$ $
\lambda, \lambda_1, \mu, \mu_1 \in F.$
\end{Definition}

 \begin{Lemma}\label {6.10}  If $0\not= I \lhd \lhd \widehat {L}(\stackrel {\circ} {L})$ and $Fc \nsubseteq I$,
 there exists $B \lhd F[t, t^{-1}]$ such that $0\not=  B \otimes \stackrel {\circ} {L} \subseteq I$.
 \end {Lemma}
 {\bf Proof.}
(i) Assume $I=I_1 \lhd I_2 \lhd \cdots \lhd I_n =
\widehat{L}(\stackrel {\circ}
 {L})$. We show this by induction on $n$. It is clear when $n=1$.
 Now $n>1$ and assume that  there exists $B_{2} \lhd F[t, t^{-1}]$ such that $0\not=
  B_2 \otimes \stackrel {\circ} {L} \subseteq I_2$.
 Let $f(t) \in B_2$ and $\partial f(t)
  >1$ (i.e. the degree of $f(t)>1$ ).

(ii)  We first show that $I \not= Fc +Fd.$ If $I = Fc +Fd$, then $d
\in I$.
  Since $[d, f(t) \otimes E_\alpha] =
  df(t) \otimes E_\alpha \not=0$
  and $ df(t) \otimes E_\alpha  \not\in I.$ This is a contradiction.

(iii) We next show that there exists $ 0\not= v \in I\cap (L
(\stackrel {\circ} {L})+Fc)$.

$(1^\circ)$ If  there exists $ 0\not= u \in I$ and $u = h(t) \otimes
x + \lambda c + \mu d$ with $0\not= h(t) \otimes x \in L (\stackrel
{\circ} {L})$, then there exists $y\in \stackrel {\circ} {L}$ such
that $[x, y] \not=0$ in $\stackrel {\circ} {L}$. Thus $0\not= [u,
h(t) \otimes y] = h(t)^2 \otimes [x, y] + dh(t) \otimes y + \psi
(h(t) \otimes x, h(t) \otimes y)c \in I\cap (L (\stackrel {\circ}
{L})+Fc) .$

$(2^\circ)$ Assume that   there exists $0\not= u  \in I$ and $u =
u_1 + \lambda c + \mu d$ with $l(u_1) >1$,   $u_1 \in L(\stackrel
{\circ} {L})$ and $\mu \not=0$.  Let $ u_1= \sum _{\alpha \in \Phi
^+} (u_\alpha \otimes e _\alpha + u_{-\alpha} \otimes e _{-\alpha} )
+ \sum _{i=1}^n u_i' \otimes H_i $ .

(a) If $v_\alpha =0$ for all $\alpha \in \Phi,$ then $0\not= [u,
f(t) \otimes H_1] \in I \cap (L (\stackrel {\circ} {L})+Fc)$.

(b) If there exists $i_0$ such that  $v_{i_0}' \not=0$ and  $v_j'
=0$ for all $j\not= i_0$, then $0\not= [u, f(t) \otimes H_{1}] \in I
\cap (L (\stackrel {\circ} {L})+Fc)$.

(c) If   $v_i'=0$   for all $i$, then there exists $\alpha _0 \in
\phi$ such that $ v_{\alpha _0 \not=0}.$ Consequently,  $0\not= [u,
f(t) \otimes e_{-\alpha _0}] \in I \cap (L (\stackrel {\circ}
{L})+Fc)$.

(iv) By the same method as the proof of Lemma \ref {6.5}, we can
show that there exists  $0\not= u = a\otimes x \in {L}(\stackrel
{\circ} {L})$  and $\lambda \in F$ such that $u+ \lambda c \in I.$
Considering Lemma \ref {6.8}, we have $0\not= u = a\otimes x \in I$.

(v) Using (ii) in  the proof of Lemma \ref {6.5}, we can complete
the proof. $\Box$

 \begin{Theorem}\label {6.11} $r_b (\widehat {L}(\stackrel {\circ} {L}))=Fc.$ \end{Theorem}
{\bf Proof.} Let $I : = r_b (\widehat {L}(\stackrel {\circ} {L}))$.
 It is clear $Fc \subseteq r_b (\widehat {L}(\stackrel {\circ} {L}))$.
 If $[\widehat {L}(\stackrel {\circ} {L}), I]=0$, then $I \subseteq Fc$ by \cite [Pro. 7.1B]{Wa02}
  and our claim holds.
If $[\widehat {L}(\stackrel {\circ} {L}), I] \not=0.$ then $I \not=
Fc$. Set $\bar L := \widehat {L}(\stackrel {\circ} {L})/ Fc.$ By
Lemma \ref {1.4}, $r_b(\bar L) \not=0$. It follows from Proposition
\ref {3.4} that there exists $0\not=A_1 \lhd A_2 \lhd \cdots \lhd
A_n \lhd \bar L$ such that $A_1$ is nilpotent. Consequently,
 there exists $0\not=B_1 \lhd B_2 \lhd \cdots \lhd B_n \lhd
 \widehat {L}(\stackrel {\circ} {L})$ such that $A_i= B_i/Fc$ for $1 \le i \le n$
  and $B_1^m\subseteq Fc$ for certain $m$. Thus $B_1$ is a non-zero  nilpotent accessible
   ideal of  $\widehat {L}(\stackrel {\circ} {L})$. By  Lemma \ref {6.10}, there exists
    $0\not= B \lhd F[t, t^{-1}]$ such that $ B\otimes \stackrel {\circ} {L} \subseteq I$.
    By the proof of Theorem \ref {6.6}, $B\otimes \stackrel {\circ} {L} $ is not nilpotent
    as a sub-set of  $ L(\stackrel {\circ} {L})$.   $B\otimes \stackrel {\circ} {L} $ is not nilpotent
    as a sub-set of  $ \widehat {L}(\stackrel {\circ} {L})$. This is a contradiction. $\Box$

\section {Appendix }\label {s7}

For every finite commutative group $G$, all results in section 1  and in section 3
 hold in the category of all $G$-colour Lie algebras (see \cite {Ka77, Sc79}) and all braided m -Lie algebras(see \cite {WZZ}).
 Therefore, all results in section 1  and in section 3  hold in the category of all super  Lie algebras (i.e. ${\mathbb Z}_2$- colour Lie algebras).

\begin{Lemma}\label {6.4}$F[t, t^{-1}]$ is a main ideal algebra (associative algebra),
 i.e. every
 ideal of $F[t, t^{-1}]$ is an ideal generated by only one element.  \end{Lemma}

 {\bf Proof.}  If $0\not= I \lhd A,$ then
 there exists $f(t) \in I \cap F[t]$ such that the degree $\partial f(t)$ of  $f(t)$ is the minimal of polynomials in $F[t]\cap I$. Consequently, $I = (f(t))$, which is the ideal generated by $f(t)$ in $F[t, t^{-1}]$. $\Box$

\end {document}